# Several Conclusions on another site setting problem

## Yuyang Zhu

*Department of Mathematics & Physics, Hefei University, Hefei* 230601, *P. R. China*
*E-mail: zhuyy@hfuu.edu.cn*

**Abstract:** Let $S = \{A_1, A_2, \cdots, A_n\}$ be a finite point set in m-dimensional Euclidean space $E^m$, and $\|A_i A_j\|$ be the distance between $A_i$ and $A_j$. Define

$$\sigma(S) = \sum_{1 \leq i < j \leq n} \|A_i A_j\|, \quad D(S) = \max_{1 \leq i < j \leq n}\{\|A_i A_j\|\}, \quad \omega(m,n) = \frac{\sigma(S)}{D(S)},$$

$$\sup \omega(m,n) = \max\left\{ \frac{\sigma(S)}{D(S)} \middle| S \subset E^m, |S| = n \right\}.$$

This paper proves that, for any point $P$ in an n-dimensional simplex $A_1 A_2 \cdots A_{n+1}$ in Euclidean space,

$$\sum_{i=1}^{n+1} \|PA_i\| \leq \sup_{i_t, j_t \in \{1,2,\cdots,n+1\}} \left\{ \sum_{t=1}^{n} \|A_{i_t} A_{j_t}\| \right\}.$$

By using this inequality and several results in differential geometry this paper also proves that

$$\sup \omega(2,4) = 4 + 2\sqrt{2 - \sqrt{3}}, \quad \sup \omega(n, n+2) \geq C_{n+1}^2 + 1 + n\sqrt{2\left(1 - \sqrt{\tfrac{n+1}{2n}}\right)}.$$

**MR(2000) Subject Classification** 52C35

## 1 Introduction

Let $S = \{A_1, A_2, \cdots, A_n\}$ be a finite point set in m-dimensional Euclidean space $E^m$, and $\|A_i A_j\|$ be the distance between $A_i$ and $A_j$. Define

$$\sigma(S) = \sum_{1 \leq i < j \leq n} \|A_i A_j\|, \quad d(S) = \min_{1 \leq i < j \leq n}\{\|A_i A_j\|\}, \quad \mu(m,n) = \frac{\sigma(S)}{d(S)},$$

$$\inf \mu(m,n) = \min\left\{ \frac{\sigma(S)}{d(S)} \middle| S \subset E^m, |S| = n \right\}.$$

It is a site setting problem to evaluate $\inf \mu(m,n)$. It has been proved with convex geometry theory that [1]

$$\inf \mu(n, n+2) = C_{n+2}^2 - 1 + 2\sqrt{\tfrac{n+1}{2n}},$$

It has also been proved that [2] $\inf \mu(2,5) = 9 + 2\sqrt{3}$. That leads to the following problem. By defining

$$D(S) = \max_{1 \leq i < j \leq n} \{\|A_i A_j\|\}, \quad \omega(m,n) = \frac{\sigma(S)}{D(S)},$$

To evaluate $\sup \omega(m,n) = \max\left\{\frac{\sigma(S)}{D(S)} \middle| S \subset E^m, |S| = n\right\}$ is a conjugate problem of the previous one. These problems enrich the discrete geometry theory and also have real world applications [1-12]. First this paper gives a geometry inequality, and then it proves that $\sup \omega(2,4) = 4 + 2\sqrt{2 - \sqrt{3}}$, $\sup \omega(n, n+2) \geq C_{n+1}^2 + 1 + n\sqrt{2\left(1 - \sqrt{\frac{n+1}{2n}}\right)}$.

## 2   Lemmas

**Lemma 2.1**  $\forall \; P \in \triangle A_1 A_2 A_3 \subset E^2$, the following inequality holds:

$$\sum_{i=1}^{3} \|PA_i\| \leq \sup_{i_t, j_t \in \{1,2,3\}} \left\{ \sum_{t=1}^{2} \|A_{i_t} A_{j_t}\| \right\}.$$

**Proof**  Let $\tau = \triangle A_1 A_2 A_3$. **(I)** If $P \in \partial \tau$, without losing generality, assume point $P$ is on edge $A_1 A_2$. Let $\angle A_1 P A_3 = \alpha_1$, $\angle A_2 P A_3 = \alpha_2$. Since $\angle A_1 + \angle A_2 < \pi = \alpha_1 + \alpha_2$, $\exists i \in \{1,2\}$, $\angle A_i < \alpha_i$ holds. Without losing generality, assume $\angle A_1 < \alpha_1$. Since a greater angle in a triangle corresponds to a longer edge, $\|PA_3\| \leq \|A_1 A_3\|$, and also considering $\|PA_1\| + \|PA_2\| = \|A_1 A_2\|$, we have

$$\sum_{i=1}^{3} \|PA_i\| \leq \|A_1 A_3\| + \|A_1 A_2\| \leq \sup_{i_t, j_t \in \{1,2,3\}} \left\{ \sum_{t=1}^{2} \|A_{i_t} A_{j_t}\| \right\}.$$

**(II)**  If point $P \in \text{int } \tau$, without losing generality, assume $\min_{i \in \{1,2,3\}}\{\angle A_i\} = \angle A_1$. Draw a line parallel to $A_1 A_2$ through $P$ and let it intersect with $A_1 A_3$ and $A_3 A_2$ on points $E$ and $E'$ respectively. Also draw a line parallel to $A_1 A_3$ and let it intersect with $A_1 A_2$ and $A_3 A_2$ on points $F$ and $F'$ respectively, and then we get

$$\|A_1 E\| + \|A_1 F\| = \|A_1 E\| + \|EP\| > \|PA_1\|$$

and

$$\angle A_3 PE = \angle PA_3 E' + \angle PE' A_3 = \angle PA_3 E' + \angle A_2 > \angle A_1 = \angle A_3 EP.$$

Thus $\|PA_3\| < \|EA_3\|$ in $\triangle A_3EP$, and similarly we can prove $\|PA_2\| < \|FA_2\|$. Therefore

$$\sum_{i=1}^{3}\|PA_i\| < \|A_1E\| + \|A_1F\| + \|EA_3\| + \|FA_2\|$$

$$= \|A_1A_3\| + \|A_1A_2\| = \sup_{i_t,j_t \in \{1,2,3\}}\left\{\sum_{t=1}^{2}\|A_{i_t}A_{j_t}\|\right\}.$$

Lemma 2.1 can be generalized with convex geometry theory and similar approaches of proving lemma 2.1, and the following conclusion can be drawn:

**Lemma 2.2** Let $\tau$ be a n-dimensional simplex in n-dimensional Euclidean space $E^n$ with vertices $A_1, A_2, \cdots, A_{n+1}$. For $\forall P \in \tau$, it holds that

$$\sum_{i=1}^{n+1}\|PA_i\| \leq \sup_{i_t,j_t \in \{1,2,\cdots,n+1\}}\left\{\sum_{t=1}^{n}\|A_{i_t}A_{j_t}\|\right\}.$$

**Proof** In the following we only prove the case in which $P \in \partial\tau$. For the case of $P \in \text{int}\,\tau$, we can prove it with similar approach seen in the proof of Lemma 2.1. If $n = 2$, according to lemma 2.1 the conclusion holds. Assume it also holds in the case of $n = k(k > 1)$. Then in the case of $n = k+1$, since $P \in \partial\tau$, there must exist a $k-1$-dimensional simplex $\tau_{k-1}$ such that $P \in \tau_{k-1}$. Without losing generality, assume $\tau_{k-1}$ = simplex $A_1A_2\cdots A_k$, using inductive assumption we get

$$\sum_{i=1}^{k}\|PA_i\| \leq \sup_{i_t,j_t \in \{1,2,\cdots,k\}}\left\{\sum_{t=1}^{k-1}\|A_{i_t}A_{j_t}\|\right\}.$$

On the other hand, since $\tau$ is a simplex and $P \in \tau$,

$$\|PA_{k+1}\| \leq \max\{\|A_{k+1}A_1\|, \|A_{k+1}A_2\|, \cdots, \|A_{k+1}A_k\|\}$$

$$= \max_{s \in \{1,2,\cdots,k\}}\{\|A_{k+1}A_s\|\}.$$

then

$$\sum_{i=1}^{k+1}\|PA_i\| = \sum_{i=1}^{k}\|PA_i\| + \|PA_{k+1}\| \leq \sup_{i_t,j_t \in \{1,2,\cdots,k\}}\left\{\sum_{t=1}^{k-1}\|A_{i_t}A_{j_t}\|\right\}$$

$$+ \max_{s \in \{1,2,\cdots,k\}}\{\|A_{k+1}A_s\|\} \leq \sup_{i_t,j_t \in \{1,2,\cdots,k+1\}}\left\{\sum_{t=1}^{k}\|A_{i_t}A_{j_t}\|\right\}.$$

Using complete induction we conclude that lemma 2.2 is proved if $P \in \partial\tau$.

**Lemma 2.3**[5],[13] If any point on a curved surface $\Sigma$ has Gauss curvature $K > 0$, then the extreme point of the surface must be a maximum point.

## 3 Theorems and proof

**Theorem 3.1** $\sup \omega(2,4) = 4 + 2\sqrt{2-\sqrt{3}}$.

**Proof** If three of four points $A_1, A_2, A_3, A_4$ in a plane are collinear, it is straightforward to prove that

$$\omega(2,4) < 5 < 4 + 2\sqrt{2-\sqrt{3}}.$$

So the only problem left is the case that no three points are collinear. If one point $A_4$ is in the simplex composed of the other three points $\Delta A_1 A_2 A_3$, then according to lemma 2.1 or lemma 2.2,

$$\sum_{i=1}^{3} \|A_4 A_i\| \leq \sup_{i_t, j_t \in \{1,2,3\}} \left\{ \sum_{t=1}^{2} \|A_{i_t} A_{j_t}\| \right\}, \text{ thus}$$

$$\omega(2,4) = \frac{\sigma(S)}{D(S)} = \left( \sum_{1 \leq i < j \leq 4} \|A_i A_j\| \right) \div \max_{1 \leq i < j \leq 4} \{\|A_i A_j\|\}$$

$$= \left( \sum_{i=1}^{3} \|A_4 A_i\| + \sum_{1 \leq i < j \leq 3} \|A_i A_j\| \right) \div \max_{1 \leq i < j \leq 4} \{\|A_i A_j\|\}$$

$$< \left( \sup_{i_t, j_t \in \{1,2,3\}} \left\{ \sum_{t=1}^{2} \|A_{i_t} A_{j_t}\| \right\} + \sum_{1 \leq i < j \leq 3} \|A_i A_j\| \right) \div \max_{1 \leq i < j \leq 4} \{\|A_i A_j\|\}$$

$$< \left( 5 \max_{1 \leq i < j \leq 4} \{\|A_i A_j\|\} \right) \div \max_{1 \leq i < j \leq 4} \{\|A_i A_j\|\} = 5 < 4 + 2\sqrt{2-\sqrt{3}}.$$

In the following we prove that $A_1, A_2, A_3, A_4$ are vertices of planar convexity, and

$$\omega(2,4) < 4 + 2\sqrt{2-\sqrt{3}}.$$

（I） If any three points of $A_1, A_2, A_3, A_4$ are vertices of a regular triangle, and the side length of the triangle $D(S) = \max_{1 \leq i < j \leq 4} \{\|A_i A_j\|\}$, then in the following we prove $\omega(2,4) < 4 + 2\sqrt{2-\sqrt{3}}$. Without losing generality, let $D(S) = 1$, vertices of the regular triangle be $A_1, A_2, A_3$. Since quadrilateral $A_1 A_2 A_3 A_4$ is convex, $A_4$ is outside the regular triangle and the distances between it and $A_1, A_2, A_3$ are all $\leq D(S) = 1$. Let $A_1, A_2, A_3$ be centers and draw arc $\overset{\frown}{A_2 A_3}$, $\overset{\frown}{A_3 A_1}$, $\overset{\frown}{A_1 A_2}$ with radius 1 respectively, then these three arcs form three bow shapes with the edges of the regular triangle, and $A_4$ must fall into one of them.

Regarding symmetry of this problem, assume $A_4$ falls into the bow shape with $A_2 A_3$ as the bowstring. Since the edge lengths of regular triangle $A_1 A_2 A_3 = D(S) = 1$, let us assume the coordinates of three vertices of the triangle are $A_1(0,0)$, $A_2(1,0)$, $A_3(\frac{1}{2}, \frac{\sqrt{3}}{2})$, and that of $A_4$ be $(x, y)$. According to the previous statements,

$$\Omega = \left\{(x, y) \big| \tfrac{1}{2} < x < 1,\ 0 < y < \tfrac{\sqrt{3}}{2},\ y + \sqrt{3}x - \sqrt{3} \geq 0,\ x^2 + y^2 \leq 1\right\},$$

which is the bow shape where $A_4$ falls into. Since $\|A_1 A_4\| \leq D(S) = 1$, the remaining task is to prove the maximum of $\|A_3 A_4\| + \|A_2 A_4\|$ is $4 + 2\sqrt{2 - \sqrt{3}}$. Let

$$f(x, y) = \|A_3 A_4\| + \|A_2 A_4\|,$$

and we have

$$f(x, y) = \sqrt{(x - \tfrac{1}{2})^2 + (y - \tfrac{\sqrt{3}}{2})^2} + \sqrt{(x - 1)^2 + y^2},$$

and

$$r = \frac{\partial^2 f}{\partial x^2} = \frac{(y - \tfrac{\sqrt{3}}{2})^2}{\left[(x - \tfrac{1}{2})^2 + (y - \tfrac{\sqrt{3}}{2})^2\right]^{\frac{3}{2}}} + \frac{y^2}{\left[(x - 1)^2 + y^2\right]^{\frac{3}{2}}},$$

$$t = \frac{\partial^2 f}{\partial y^2} = \frac{(x - \tfrac{1}{2})^2}{\left[(x - \tfrac{1}{2})^2 + (y - \tfrac{\sqrt{3}}{2})^2\right]^{\frac{3}{2}}} + \frac{(x - 1)^2}{\left[(x - 1)^2 + y^2\right]^{\frac{3}{2}}},$$

$$s = \frac{\partial^2 f}{\partial x \partial y} = -\left(\frac{(x - \tfrac{1}{2})(y - \tfrac{\sqrt{3}}{2})}{\left[(x - \tfrac{1}{2})^2 + (y - \tfrac{\sqrt{3}}{2})^2\right]^{\frac{3}{2}}} + \frac{(x - 1)y}{\left[(x - 1)^2 + y^2\right]^{\frac{3}{2}}}\right).$$

Since $A_4$ is not collinear with $A_2 A_3$, or equivalently $\dfrac{y - \tfrac{\sqrt{3}}{2}}{x - \tfrac{1}{2}} \neq \dfrac{y}{x - 1}$, according to Cauchy inequality, $rt - s^2 > 0$. According to differential geometry theory on curved surface, the function $f(x, y)$ has its Gauss curvature $K > 0$ everywhere in $\Omega$. According to lemma2.3, the extreme point of $f(x, y)$ in $\Omega$ must be a maximum point. In the following we evaluate the stable point of $f(x, y)$ in $\Omega$.

$$\begin{cases} \dfrac{\partial f}{\partial x} = \dfrac{x-\frac{1}{2}}{\left[(x-\frac{1}{2})^2+(y-\frac{\sqrt{3}}{2})^2\right]^{\frac{1}{2}}} + \dfrac{x-1}{\left[(x-1)^2+y^2\right]^{\frac{1}{2}}} = 0, \\ \dfrac{\partial f}{\partial y} = \dfrac{y-\frac{\sqrt{3}}{2}}{\left[(x-\frac{1}{2})^2+(y-\frac{\sqrt{3}}{2})^2\right]^{\frac{1}{2}}} + \dfrac{y}{\left[(x-1)^2+y^2\right]^{\frac{1}{2}}} = 0. \end{cases} \quad (1)$$

Solve (1) and we get

$$(1-x)(y-\tfrac{\sqrt{3}}{2}) = (\tfrac{1}{2}-x)y, \qquad (2)$$

or

$$(1-x)(y-\tfrac{\sqrt{3}}{2}) = (x-\tfrac{1}{2})y. \qquad (3)$$

Since $A_4 \in \Omega$, the left hand side of (3) $<0$ while the right hand side $>0$, then (3) does not hold. Only (2) holds and according to (2),

$$y + \sqrt{3}x - \sqrt{3} = 0,$$

which is the equation of line passing through points $A_2$ and $A_3$. Since $A_4$ is not collinear with $A_2 A_3$, function $f(x,y)$ has no stable point in $\Omega$. According to the continuity of $f(x,y)$ in $\Omega$, $f(x,y)$ must reach its maximum on the boundary of $\Omega$, or $\partial \Omega$. According to previous statement, if $f(x,y)$ reaches its extreme value on $\partial\Omega$, then we can evaluate its maximum there. Since

$$\partial\Omega = \left\{(x,y)\,\middle|\, y+\sqrt{3}x-\sqrt{3}=0, \tfrac{1}{2}\le x\le 1\right\} \cup$$
$$\left\{(x,y)\,\middle|\, x^2+y^2=1, \tfrac{1}{2}\le x\le 1,\ 0\le y\le \tfrac{\sqrt{3}}{2}\right\},$$

it is obvious that if point $A_4(x,y)$ is collinear with $A_2 A_3$, then $y+\sqrt{3}x-\sqrt{3}=0$, and if $\tfrac{1}{2}\le x\le 1$, then $f(x,y)=1<2\sqrt{2-\sqrt{3}}$.

If $A_4(x,y) \in \left\{(x,y)\,\middle|\, x^2+y^2=1, \tfrac{1}{2}\le x\le 1,\ 0\le y\le \tfrac{\sqrt{3}}{2}\right\}$, then according to Lagrange multiplier method the only maximum point is $(\tfrac{1}{2},\tfrac{\sqrt{3}}{2})$, where the maximum of $f(x,y)$ is reached, or $\max\limits_{(x,y)\in\partial\Omega}\{f(x,y)\} = f(\tfrac{1}{2},\tfrac{\sqrt{3}}{2}) = 2\sqrt{2-\sqrt{3}}$. Hence, we proved that in the convex quadrilateral if there are three vertices which are also the vertices of a regular triangle and let $D(S)$ equals to the side length of the triangle, then

$$\omega(2,4) = \sum_{1 \le i < j \le 3} \|A_i A_j\| + \|A_1 A_4\| + f(x,y) \le 3 + 1 + 2\sqrt{2-\sqrt{3}} = 4 + 2\sqrt{2-\sqrt{3}}.$$

The equality holds if and only if $A_4$ $(x, y)$ bisects the arc.

(II) Since there must be two points between which the distance reaches maximum in a convex quadrilateral, without losing generality, let $\|A_1 A_2\| = D(S)$ and the coordinates of $A_1$ and $A_2$ be $A_1(0,0)$ and $A_2(1,0)$ respectively, while other two points be in the region

$$\Omega_1 = \{(x,y) | x^2 + y^2 \le 1\} \cap \{(x,y) | (x-1)^2 + y^2 \le 1\}.$$

Since any point $P$ not collinear with $A_1 A_2$ must be on a circular arc $\widehat{A_1 A_2}$. In other words, $\forall P(x,y) \in \Omega_1$ and $y \ne 0$, $\exists b \in R$, $r \in R^+$ such that $P(x,y)$ satisfies the following equation:

$$(x - \tfrac{1}{2})^2 + (y - b)^2 = r^2. \quad (4)$$

This arc is denoted by $S^1(b,r)$. In the following we prove that if $P(x,y) \in S^1(b,r)$, to make $\|A_1 P\| + \|A_2 P\|$ reaches its maximum on the circle arc, $P$ must bisect $S^1(b,r)$ and its coordinate satisfies (4), and

$$\left(\|A_1 P\| + \|A_2 P\|\right)^2 \le 2\left(\|A_1 P\|^2 + \|A_2 P\|^2\right) = 2\left(x^2 + y^2 + (x-1)^2 + y^2\right).$$

Let $L = 2\left(\|A_1 P\|^2 + \|A_2 P\|^2\right) = 2\left(x^2 + y^2 + (x-1)^2 + y^2\right)$, then according to (4)

$$L = 2\left(x^2 + y^2 + (x-1)^2 + y^2\right) = 2\left(x^2 + (x-1)^2 + 2\left(b + \sqrt{r^2 + (x-\tfrac{1}{2})^2}\right)\right).$$

If

$$L'_x = 2\left[4x - 2 - 4(x - \tfrac{1}{2}) - \left(\frac{4b(x - \tfrac{1}{2})}{\sqrt{r^2 + (x-\tfrac{1}{2})^2}}\right)\right] = 0,$$

since $P \in \Omega_1$, $b \ne 0$ and $x = \tfrac{1}{2}$, in other words, if $P$ bisects $S^1(b,r)$, then $L$ reaches its maximum. If $x = \tfrac{1}{2}$, then $\|A_1 P\| = \|A_2 P\|$, hence $\|A_1 P\| + \|A_2 P\| = \sqrt{L}$ also reaches its maximum. According to this conclusion, only if points $A_3, A_4$ are on the midnormal of $A_1 A_2$, then $\left(\|A_1 A_3\| + \|A_2 A_3\|\right) + \left(\|A_1 A_4\| + \|A_2 A_4\|\right)$ can reach its maximum. On the other hand, since $\|A_3 A_4\| \le D(S) = 1$, only if $\|A_3 A_4\| = D(S) = 1$ and points $A_3, A_4$ are on the midnormal of $A_1A_2$, then $\sum_{1\le i<j\le 4}\|A_iA_j\|$ can reaches its maximum. Thus we can let the coordinates of $A_3, A_4$ be $A_3(\frac{1}{2}, y_1)$ and $A_4(\frac{1}{2}, y_2)$ respectively, where $|y_1 - y_2| = 1$ and $y_1, y_2 \in (-\frac{\sqrt{3}}{2}, \frac{\sqrt{3}}{2})$, and use similar approach as in (I) to prove that $\omega(2,4) \le 4 + 2\sqrt{2-\sqrt{3}}$.

To conclude, $\sup \omega(2,4) \le 4 + 2\sqrt{2-\sqrt{3}}$, and according to (I) the equality does conditionally hold, hence

$$\sup \omega(2,4) = 4 + 2\sqrt{2-\sqrt{3}}.$$

The theorem is proved.

From definition we can also conclude that

**Theorem 3.2** $\sup \omega(n, n+1) = C_{n+1}^2$.

Using the method of constitution, the following theorem can be proved:

**Theorem 3.3** $\sup \omega(n, n+2) \ge C_{n+1}^2 + 1 + n\sqrt{2\left(1 - \sqrt{\frac{n+1}{2n}}\right)}$.

**Proof** Consider $n+2$ points $A_1, A_2, \cdots, A_{n+1}, A_{n+2}$ in n-dimensional Euclidean space $E^n$. If $A_1, A_2, \cdots, A_{n+1}$ are the vertices of n-dimensional regular simplex with unit side length, then draw a n-dimensional hyper-sphere $O$ whose center is $A_1$ and with unit radius length. $n-1$-dimensional simplex $A_2 A_3 \cdots A_{n+1}$ gives a $n-1$-dimensional hyper plane and divide $O$ into two parts of sphere with different sizes. Denote the median point of the smaller one face by $A_{n+2}$, then according to convex geometry theory we have

$$\|A_iA_j\| = 1 (\forall i, j \in \{1, 2, \cdots, n+1\}, i \ne j), \quad \|A_1A_{n+2}\| = 1,$$

$$\|A_iA_{n+2}\| = \sqrt{\left(1 - \sqrt{\frac{n+1}{2n}}\right)^2 + \left(\frac{n-1}{n}\sqrt{\frac{n}{2(n-1)}}\right)^2}$$

$$= \sqrt{2\left(1 - \sqrt{\frac{n+1}{2n}}\right)} \ (i = 2, 3, \cdots, n+1),$$

and $\sqrt{2\left(1 - \sqrt{\frac{n+1}{2n}}\right)} < 1$, we conclude that

$$\omega(n, n+2) = \sum_{1\le i<j\le n+1}\|A_iA_j\| + \|A_1A_{n+2}\| + \sum_{i=2}^{n+1}\|A_iA_{n+2}\|$$

$$= C_{n+1}^2 + 1 + n\sqrt{2\left(1 - \sqrt{\frac{n+1}{2n}}\right)},$$

thus

$$\sup \omega(n, n+2) \geq C_{n+1}^2 + 1 + n\sqrt{2\left(1 - \sqrt{\tfrac{n+1}{2n}}\right)}.$$

The theorem is proved.

## 4    Conclusions and guesses

The proof of $\sup \omega(2,4) = 4 + 2\sqrt{2-\sqrt{3}}$ is more complicated compared with the proof of $\inf \mu(2,4) = 5 + \sqrt{3}$ ([9]). Therefore the evaluation of $\sup \omega(m,n)$ may be more difficult than that of $\inf \mu(m,n)$. According to the conclusion of [10], if $n(n > m)$ points are vertices of a convexity in $m$-dimensional Euclidean space, and their distribution is such that $\inf \mu(m,n)$ can be reached, then the distances between neighboring points must be equal. However, this is not case for $\sup \omega(m,n)$, which is shown in the proof of $\sup \omega(2,4) = 4 + 2\sqrt{2-\sqrt{3}}$.

According to theorem 2.2 and [1], $\sup \omega(n,n+1) = \inf \mu(n,n+1) = C_{n+1}^2$. On the other hand according to the conclusion in [9] $\inf \mu(2,4) = 5 + \sqrt{3} > \sup \omega(2,4)$, can prove $\inf \mu(2,5) > \sup \omega(2,5)$.

**Guess 4.1**   If $n > m+1$, and $m \geq 2$, then $\inf \mu(m,n) > \sup \omega(m,n)$.

**Guess 4.2**   $\sup \omega(n, n+2) = C_{n+1}^2 + 1 + n\sqrt{2\left(1 - \sqrt{\tfrac{n+1}{2n}}\right)}.$

**References**
[1] YuYang ZHU. Several Results on a Problem of Setting Site. *Acta Mathematica Sinica, Chinese Series A,* 2011, **54**(4):669-676.
[2] Zhu Yuyang. The Proof of the Conjecture on a Problem of Setting Site. *COMM.ON APPL. MATH. AND COMPUT.* 2010, Vol.24 No.2:101-112(in Chinese).
[3] Peter Brass, William Moser, Jănos Pach. Research Problems in Discrete Geometry [M]. 2005 Springer Science +Business Inc. 2－159, 443－450.
[4] J. Pach, P. Agarwal. Combinatorial Geometry. *John wiley ＆ Sons, Ine.1995.*
[5] ZHU Yuyang. A Problem of Site Setting.*Journal of University of Science And Technology of China,* 2011, Vol.41 No.6:480-491.
[6] Hong Y, Wang G Q, Tao Z S. On a Problem of Heilbronn Type in Higher Dimensional Space [J]. *Acta Math. Sinica , Chinese Series A*, 1997,40(1): 144－153.
[7] Tao Z S, Hong Y. On a Problem of Heilbronn`s Type in $R^3$ [J]. *Acta Math. Sinica , Chinese Series A*,2000,43(5):797-806.
[8] Yang L, Zhang J Z, Zeng Z B. On the Heilbronn numbers of triangular regions[J]. *Acta Math. Sinica , Chinese Series A*,1994,37(5):678-689.
[9] Zhu Yu-yang. Discussion on the Problem of Setting up Sites[J]. *J. Hefei institute Education*, 2002, **19**(2): 1-4.

[10] Zhu Yu-yang. Some Answers to a Question on Discrete Point-Set Extremum. *J. Hefei university*, 2004, **21**(4):1-4.

[11] Zhu Yu-yang, Zhang Xia, Chu Zhao-hui. An Extreme-Value Problem of Plane Point Set. *J. Hefei university*, 2006, Vol.16 No.2:1-4.

[12] Zhu Yu-yang. Introduction to Discrete and Combinatorial Geometry. University of Science And Technology of China Press, Hefei, 2008.

[13] Peng Jia-gui, Chen Qing. Differential Geometry. Higher Education Press, Beijing, 2002.